\documentclass[11pt]{pjm}

\usepackage{amsmath,amsthm,amsfonts,latexsym,amscd}

\usepackage{eucal}

\usepackage{epsfig,graphicx}

\swapnumbers

\theoremstyle{plain}
    \newtheorem{theorem}{Theorem}[section]
    \newtheorem{lemma}[theorem]{Lemma}
    \newtheorem{corollary}[theorem]{Corollary}

\theoremstyle{definition}
    \newtheorem{definition}[theorem]{Definition}
    \newtheorem{example}[theorem]{Example}
    \newtheorem{notation}[theorem]{Notation}

\theoremstyle{remark}
\newtheorem{remark}[theorem]{Remark}
\newtheorem{remarks}[theorem]{Remarks}

\numberwithin{equation}{section}

\addtocounter{section}{-1}

\newcommand{\script}[1]{{\mathcal{#1}}}

\newcommand{\gs}{\sigma}
\newcommand{\ga}{\alpha}

\newcommand{\gl}{\lambda}

\newcommand{\eps}{\epsilon}
\newcommand{\up}{{\uparrow}}
\newcommand{\xv}{\mathbf x}
\newcommand{\yv}{{\mathbf y}}
\newcommand{\tv}{\mathbf t}

\newcommand{\zv}{\mathbf z}
\newcommand{\psp}{p_s^+}
\newcommand{\psm}{p_s^-}
\newcommand{\pspm}{p_s^\pm}
\newcommand{\ptspp}[1]{p_{s_{#1},\tv_{#1}}^+}
\newcommand{\ptsmm}[1]{p_{s_{#1},\tv_{#1}}^-}
\newcommand{\qtspp}[1]{q_{s_{#1},\tv_{#1}}^+}
\newcommand{\qtsmm}[1]{q_{s_{#1},\tv_{#1}}^-}
\newcommand{\ptsppmm}[1]{p_{s_{#1},\tv_{#1}}^\pm}
\newcommand{\qtsppmm}[1]{q_{s_{#1},\tv_{#1}}^\pm}

\newcommand{\bt}[1]{b_{\tv_{#1}}}
\newcommand{\ptsp}{p_{s,\tv}^+}

\newcommand{\ptsm}{p_{s,\tv}^-}

\newcommand{\ptspm}{p_{s,\tv}^{\pm}}

\newcommand{\pts}{p_{s,\tv}}

\newcommand{\rn}[1]{$\mathbb R^{#1}$}
\newcommand{\rnn}[1]{\mathbb R^{#1}}

\newcommand{\e}{\operatorname{ext}}

\newcommand{\conv}{\operatorname{conv}}

\newcommand{\down}{\downarrow}
\newcommand{\lan}{\langle}
\newcommand{\ran}{\rangle}
\newcommand{\spe}{$(s,\tv)$ }
\newcommand{\spk}[1]{(s_{#1},\tv_{#1})}

\renewcommand{\labelenumi}{$(\arabic{enumi})$}

\DeclareMathOperator{\spn}{span}

\begin{document}

\title[A Geometric spectral theory]{A Geometric spectral theory for
$n$-tuples of
self-adjoint operators in finite von Neumann algebras:  II}

\author[Akemann]{Charles A. Akemann}
\address{Department of Mathematics\\
University of California\\
Santa Barbara, CA 93106}
\email{\tt akemann@math.ucsb.edu}
\urladdr{http://www.math.ucsb.edu/\textasciitilde akemann/}

\author{Joel Anderson}
\address{Department of Mathematics\\
Pennsylvania State University\\
University Park, PA 1680}
\email{\tt anderson@math.psu.edu}


\thanks{Each author was supported in part by The National Science
Foundation during the
preparation of this paper.}







\date{Received date / Revised version date}

\begin{abstract}
Given an $n$-tuple $\{b_1, ..., b_n\}$ of self-adjoint
operators in a finite von Neumann algebra  $M$  and a faithful, normal
tracial state $\tau$ on  $M$, we define a map
$\Psi$ from   $M$ to \rn{n} by $\Psi(a) = (\tau(a),\tau(b_1a), \dots,
\tau(b_na))$.  The image of the positive part of the unit ball under $\Psi$
is called
the {\bf spectral scale} of $\{b_1, .., b_n\}$ relative to
$\tau$ and is denoted by $B$.   In a previous paper with Nik Weaver we
showed that the geometry
of
$B$ reflects spectral data for real linear combinations of the operators
\{$b_1, .., b_n$\}. For example, we showed that an exposed face in $B$ is
determined by
a certain pair of spectral projections of a real linear combination of the
$b_i$'s.
In the present paper we extend this study to faces that are not exposed.
In order to do
this we need to introduce a recursive method for describing faces of
compact convex sets
in \rn{n}. Using this new method, we completely describe the  structure of
arbitrary faces of $B$ in terms of $\{b_1, .., b_n\}$ and $\tau$. We also
study faces
of convex, compact sets that are exposed by more than one hyperplane of
support (we call
these {\bf sharp faces}).  When such faces appear on $B$,
they signal the existence of commutativity among linear combinations of the
operators
\{$b_1, .., b_n$\}. Although many of the conclusions of this study involve
too much
notation to fit nicely in an abstract, there are two results that give
their flavor very
well. Theorem 6.1:  If the set of extreme points of $B$ is countable, then $N =
\{b_1,\dots,b_n\}^{\prime\prime}$ is abelian.  Corollary 5.6:  $B$ has a
finite number of
extreme points if and only if
$N$ is abelian and has finite dimension.
\end{abstract}

\maketitle

\section{Introduction}\label{intro}

In this work we continue the investigation begun in \cite{AAW} of the spectral
relationships within an $n$-tuple of self--adjoint operators in a finite von
Neumann algebra as reflected in  the geometry of the  corresponding
spectral scale.
Let us begin by introducing some notation and  reviewing  some of the
results in
\cite{AAW}.

{\bf The following notation will apply throughout the rest of the paper without
explicit reference.}

\begin{notation}Let $M$ denote a finite von Neumann algebra equipped with a
faithful
normal tracial state $\tau$, let  $n$ denote a positive integer, and let
$b_1,\dots,b_n$
denote an $n$--tuple of self--adjoint operators in $M$. We define a map
$\Psi$ from $M$
to $\mathbb R^{n+1}$ by  the formula
\[
\Psi(a) = (\tau(a),\tau(b_1a), \dots, \tau(b_na))
\]
and write $B =\Psi(M_1^+)$, where $M_1^+ = \{a\in M: 0 \le a \le 1\}$.
Since $\tau$ is
normal and $M_1^+$ is weak* compact and convex, it follows that $B$ is a
compact, convex
subset of $\mathbb R^{n+1}$. We call  $B$ the {\bf spectral scale} of the
$b_i$'s {\bf
relative to $\tau$}.

We use parentheses [$(,)$] to denote vectors in \rn{n} and angle brackets
[$\lan, \ran$]
to denote inner products.

\medskip
\begin{enumerate}

\item We let $N$ denote the von Neumann subalgebra of  $M$ generated by
$b_1,\dots,b_n$ and the identity.

\item If $p$ and $q$ are projections in $M$ and $p \le q$, then the {\bf
order interval}
that they determine is
\[
[p,q] = \{a\in M^+_1: p \le a \le q\}.
\]

\item For each nonzero vector $\tv = (t_1, \dots, t_n) \in \mathbb R^n$ we
write
\[
b_{\tv} = t_1b_1 + \mathbb \cdots + t_nb_n.
\]

\item  By a {\bf spectral pair} we mean a pair of the form $(s,\tv)$, where
$s$ is a real number and $\tv$ is a nonzero vector in \rn{n}.  Our standard
way of embedding $\mathbb R^{n}$ into $\mathbb R^{n+1}$ is $\tv \rightarrow
(0,\tv)$.
\item If $(s,\tv)$ is a spectral pair, then $p_{s,\tv}^+$ and $\ptsm$
denote the spectral
projections of $b_{\tv}$ determined by the intervals $(-\infty, s]$ and
$(-\infty, s)$. We call $\ptspm$ the {\bf spectral interval projections}
determined by $s$
and $\bt{}$.
\item  If $(s,\tv)$ is a spectral pair and $\ga$ is a real number, then
$P(\tv, s,\alpha)$ denotes  the hyperplane in $\mathbb R^{n+1}$ defined by
the formula
\[
-sx_0 + t_1x_1 + \cdots + t_nx_n = \alpha
\]
and we write $P^\up(s,\tv,\alpha)$for the half--space defined by the inequality
\[
-sx_0 + t_1x_1 + \cdots + t_nx_n \geq \alpha.
\]
\item If $0 \le s \le 1$, then the {\bf isotrace slice} of $B$ at $s$ is by
definition
\[
I_s = \{\xv = (x_0,\dots,x_n)\in B: x_0 = s\}
\]
\end{enumerate}
\end{notation}
The notion of an isotrace slice may seem artificial, but, as we shall show in a
forthcoming paper, for $n=2$ it essentially contains the more familiar
concept of the
numerical range of the operator $c = b_1+ib_2$.  Isotrace slices are
particularly
useful when there are two operators so that $B \subset \rnn{3}$ as will be
seen in
\S 6 of this paper.

As shown in \cite{AAW}, the geometry (i.e. the facial structure) of the
spectral scale $B$ contains information about the spectrum and spectral
projections of
real linear combinations of the $b_1,\dots,b_n$.  In this paper we complete
this part of
the theory by giving an exact description of an arbitrary face of $B$ in
terms of the
$b_i$'s.  In \cite{AAW} this was done only for the case $n =1$, and it is
useful to begin by presenting a description of these results.
\bigskip

Assume $n=1$ and $b_1= b$.  Thus, $B$ is a convex, compact subset of the
plane so that a
face in $B$ is either an extreme point or a line segment. Also, for each
real $s$ write
$p_s^+$, (resp., $p_s^-$) for the spectral projection of $b$ corresponding
to the
interval $(-\infty,s]$ (resp., $(-\infty,s)$).  Finally write $\gs(b)$ for
the spectrum
of $b$. The results obtained in \cite[Theorems 1.5 and 1.6]{AAW}  may be
summarized as
follows.

\renewcommand{\makelabel}{(1)}
\begin{enumerate}

\item $B$ lies between the lines $x = 0$ and $x = 1$ and has ``sharp
points'' at $\Psi(0) = (0,0)$ and $\Psi(1) = (1,\tau(b))$.  Thus, the
boundary of $B$ is
divided into an upper and lower boundary.

\item  The spectrum of $b$ is exactly the set of slopes of the tangent
lines to the lower
(or upper) boundary of $B$.

\item  A point $\xv$ in $B$ is an extreme point on the lower boundary  if
and only if  $\xv$ has the form
$\Psi(p_s^\pm)$ for some $s\in  \gs(b)$.  Extreme points on the upper boundary
have the form $\Psi(1 - \ptspm)$.

\item  The line segments in the lower (or upper) boundary of $B$ are in
one-to-one correspondence with the eigenvalues of $b$.  If  $s$  is an
eigenvalue for $b$
so that  $\psm < \psp$ and $F$ is the line segment on the lower boundary of
$B$ that it
determines, then $F =\Psi([\psm,\psp])$.

\item The corners in the lower (or upper) boundary of the spectral scale are in
one-to-one correspondence with gaps in the spectrum of the operator $b$.
\end{enumerate}

Let us now review the basic facts obtained  in \cite{AAW} for the higher
dimensional
case.  Since we shall refer to these facts repeatedly, it is convenient to
give them local numbers.

The following is a restatement of the results in Theorems 2.3 and 2.4 in
\cite{AAW}.
\begin{theorem} The following statements hold.
\begin{enumerate}

\item  If $\xv$ is an extreme point of $B$, then there is a projection $p$ in
$N$ such that
\[
\Psi(p) = \xv \text{ and } \Psi^{-1}(\xv)\cap M_1^+ = \{p\}.
\]
Further, $\Psi(p_{s,\tv}^+)$  and $\Psi(p_{s,\tv}^-)$ are extreme points of
$B$ for every spectral pair $(s,\tv)$.

\item We have  $\tau((b_{\tv} - s1)p_{s,\tv}^+) = \tau((b_{\tv} -
s1)p_{s,\tv}^-)$. The
hyperplane $P(s,\tv,\ga)$ is a hyperplane of support for $B$ with $B\subset
P^\up(s,\tv,\ga)$  if and only if
\[
\ga = \tau((b_{\tv} - s1)p_{s,\tv}^\pm).
\]
In this case we have  $\Psi(\ptspm) \in P(s,\tv,\ga)$.

\item  If  $\ga = \tau((b_{\tv} - s1)p_{s,\tv}^\pm)$, then   $F =
P(s,\tv,\ga)\cap B$ is
a face in $B$. Moreover,
\[
\Psi^{-1}(P(s,\tv,\alpha)) \cap M_1^+ = [p_{s,\tv}^-, p_{s,\tv}^+]
\]
and
\[
F = \Psi([p_{\tv, s}^-, p_{\tv, s}^+]).
\]
\end{enumerate}
\end{theorem}

Recall that a face $F$ of a convex set $C$ is said to be an {\bf exposed face} if there
is a hyperplane of support $P$ for $C$ such that $F = C\cap P$. If the exposed face $F$
is a single point $\xv$, and so necessarily an extreme point of $B$, we call $\xv$ an
{\bf exposed point}.  (See \S2 for a more detailed review of these notions).

Although it was not specifically noted in \cite{AAW}, Theorem 0.2  allows the complete
classification of the exposed faces in the spectral scale which we record below.

\begin{corollary}The following statements hold.
\begin{enumerate}

\item  $F$ is an exposed face in $B$ if and only if there is a spectral pair \spe such
that
\[
F = \Psi([\ptsm,\ptsp]).
\]

\item   If $\xv$ is an extreme point in the spectral scale, then $\xv$ is an exposed
point if and only if  there is a spectral pair \spe such that $\ptsm  = \ptsp$ and $\xv =
\Psi(\ptspm)$.

\item  If $F$ is a face that is not an extreme point, then $F$ is an exposed face if and
only if there is a spectral pair \spe such that $s$ is an eigenvalue for $\bt{}$.
Thus, in this case $(\ptsp-\ptsm) \neq 0$ and $\bt{}(\ptsp-\ptsm) =s(\ptsp-\ptsm)$.
\end{enumerate}
\end{corollary}

\begin{proof} The first conclusion follows immediately from parts $(2)$ and $(3)$ of
Theorem 0.2.

Suppose $\xv$ is an exposed point so that by Theorem 0.2 and part 1 of this corollary we
have
\[
\{\xv\} = B\cap P(s,\tv,\ga) =  \Psi([\ptsm,\ptsp]) = \Psi(\ptspm),
\]
where \spe is a spectral pair. Write $r = \ptsp-\ptsm$ and
\[
M_r = rMr|_r.
\]
If $r \ne 0$, then there is a nontrivial projection $r^\prime$ in $M_r$ such that
$0 < \tau(r^\prime) < 1$.  But in this case, we have
\[
\tau(\ptsm) < \tau(\ptsm + r^\prime) < \tau(\ptsp)
\]
so that  $\Psi(\ptsm) \ne \Psi(\ptsp)$ and $\Psi([\ptsm,\pts]) \ne \{\xv\}$.  As this is
a contradiction, we get $\ptsm = \ptsp$, as desired.

Conversely, if \spe is a spectral pair such that $\ptsm = \ptsp$ and $\xv =
\Psi(\ptspm)$, then
\[
B\cap P(s,\tv,\ga) = \Psi(\pspm) = \Psi([\ptsm,\ptsp])
\]
is a singleton  and so this face is an exposed point.

Finally,  by part 1 and the fact that  $\Psi$ is faithful, we have that
$(\ptsp-\ptsm)\neq 0$  exactly when the face  $F = \Psi([\ptsm,\ptsp])$ has
positive
dimension.

\end{proof}

\bigskip

 The results to be presented below concern the following topics.

{\smallskip \narrower\narrower
\noindent  $\bullet$ A general analysis of a face of a compact, convex set
in \rn{n}.  (\S2 and \S5).

\noindent  $\bullet$ A complete description of
the facial structure of the spectral scale in the general case. (\S3).

\noindent $\bullet$ An analysis of the ``corners'' of the spectral scale in
higher dimensions.
(\S4).

\noindent $\bullet$ Applications of the results in \S4 to show how geometric
properties of $B$ imply the existence of central projections in $N$. (\S5
and \S6).

\noindent
\medskip}

Let us now describe our results in more detail.   Although we defined
exposed faces of
convex compact subsets  in \rn{n}  using hyperplanes of support, there is
also an
equivalent formulation in terms of linear functionals.  From this point of
view,
a face $F$ in the (convex, compact) set $C$ is said to be {\bf exposed} if
there is
linear functional $f$  and a scalar $t$ such that $f(\xv) = t$ for each
$\xv$ in $F$ and
$f(\xv) < t$ for each  $\xv$ in  $C \setminus F$.  In this case, we say
that $f$ {\bf
exposes} $F$. Since the former definition is more geometric it fits better
with the
emphasis of this paper, but it will sometimes be convenient to use the
latter, more
algebraic, definition.

Observe that not every face need be exposed.  For example, it is easy to
construct a
convex set in two dimensions that has a face in its boundary of dimension
one such that
the end points of this line segment (which are extreme points) are not
exposed.  We wish
to make distinctions among exposed faces as follows.  The exposed face $F$
is said to be
have {\bf degree}  ${\mathbf k}$  if $k+1$  is the cardinality of the
largest linearly
independent set of linear functionals that expose  $F$.  In terms of
hyperplanes, this
means that there are $k+1$ hyperplanes of support for $F$ such that their
normal vectors
are linearly independent.  We study such faces in \S4.   

In \S3 we show that a general face of the spectral scale provides spectral
information for
linear combinations of the defining $b_i$'s and certain cut--downs of these
operators. If
$F$ is not exposed, then a condition similar to that of Corollary 0.1 (1)
holds, except
that it is first necessary to cut down by a certain spectral projection.
This analysis
allows us to complete the  characterization of the extreme points of the
spectral scale
which was begun in \cite{AAW}.

Generalized ``corners'' are studied in \S4.  A corner of a convex  planar
set is a point
on the boundary that admits more than one tangent line of support.  The natural
generalization of this notion in higher dimensions is a face that admits
more than one
hyperplane of support, i.e., a face of degree $k$ for some $1 <  k \le n $.
We call such
faces {\bf sharp faces} since they generalize sharp (i.e.,
nondifferentiable) corners on
the boundary of a 2 dimensional convex set.

Our main result on this topic can be described as follows.  In two
dimensions, a corner
of a convex set admits precisely two tangent lines of support with linearly
independent
normal vectors.  In higher dimensions, there can be much wider variation
and a new
phenomenon occurs.

Specifically, in the case of a spectral scale formed by a single operator,
a tangent
line of support for a face of the spectral scale is determined by real numbers
$s$ and $\ga$, where $s\in \gs(b)$.   In the general case where there are
$n$ $b_i$'s, the
hyperplanes of the support for a sharp face are determined by a sequence
$\spk{1},\dots,\spk{k}$ of spectral pairs and real numbers $\ga_1,\dots,\ga_k$.

For example, suppose $F$ is a sharp face which is contained in hyperplanes
of the form
$P(s_1,\tv_1,\ga_1)$ and $P(s_2,\tv_2,\ga_2)$, where the spectral pairs
$(s_1,\tv_1)$
and $(s_2,\tv_2)$ are linearly independent. If $\tv_1 =\tv_2 = \tv$ and
$s_1 < s_2$,
then, just as in two dimensions, the interval $(s_1,s_2)$ lies in
a gap in the spectrum of $\bt{}$.  On the other hand if $\{\tv_1, \tv_2\}$
is a linearly
independent set, then something new occurs. In this case, there is a
projection $r$ that
commutes with $\bt{1}$ and
$\bt{2}$ and such that
\[
\bt{1}r = s_1r \text{ and } \bt{2}r = s_2r.
\]
If there are $k$ hyperplanes containing $F$, with linearly independent vectors
$\tv_1,\dots,\tv_k$, then the projection $r$ commutes with each $\bt{i}$
and the
compression of each $\bt{i}$ to $r$ is a scalar.

Next, we observe in \S5 that if $k = n$ for the sharp face  $F$ of the previous
paragraph, then the projection $r$ commutes with each $b_i$ and so $r$ is
a nontrivial
central projection in $N$.  Put colloquially, if we can ``wobble'' a sharp
face in all $n$
``$\tv$--directions'', then the center of  $N$ is not  trivial.  \S6 is
devoted to
showing that if $B$ has a countable number of extreme points, then $N$ is
abelian.

\bigskip

\section{Examples}

\bigskip
As noted in the introduction, one can ``see'' if $N$ is abelian from the
shape of the
spectral scale in the finite dimensional case, but the situation in
infinite dimensions is
much more complicated.  In this section we present several examples that
illustrate this
point. Other examples of the spectral scale in \rn{2} were presented in
\cite{AAW} at the
end of \S1.  Examples for which  $N$ is finite dimensional and
noncommutative appear in a forthcoming paper.  We begin with an example
where $N$ is
abelian and $B$ has a countable number of extreme points.
\bigskip

\example Let $b$  denote the diagonal operator on $\ell^2(\mathbb N)$ defined
by $b(\{r_k\}) = \{r_k/k\}$, and let  $M$ be the von Neumann algebra
generated by  $b$
and the identity (so that $M = \ell^\infty(\mathbb N)$).  Also, define a
trace $\tau$ on
$M$  by $\tau (\{t_k\}) = \sum t_k/2^k$.  Since the spectrum of  $b$  is
exactly $\{1/k\}$
together with 0, it is countable.  By Theorem 2.2 the spectral scale of
$\{b,\tau\}$ is
a 2 dimensional, convex, compact set with upper boundary curve consisting
of line
segments with decreasing slopes --- specifically the slopes $\{1/k\}$. Thus
the algebra
$N$  is abelian (it equals $M$) and infinite dimensional, and the set of
extreme points
of $B$ is countable.
\endexample

\bigskip

As the next example shows, even if $N$ is abelian and all linear
combinations of $\{b_1,
b_2\}$ have countable spectrum, $B$ may still have uncountably many extreme
points.

\example  Let  $M$ and $\tau$ be as in Example 1.1 and select
  positive  real sequences  $\{s_k\}$ and $\{t_k\}$, each of which
converges to 0, and
such that the sequence $\{t_k/s_k\}$  takes on every rational value between
0  and
infinity. Also, let
$b_1$ and $b_2$ denote the corresponding diagonal operators in
$\ell^\infty(\mathbb N)$.
Since the two sequences converge to  0,  every linear combination of $b_1$
and $b_2$ is
compact.  In particular, each linear combination of  $b_1$ and $b_2$ has
countable
spectrum.

To see that $B$ has an uncountable number of extreme points, fix an
irrational number $\ga
> 0$, write $\tv = (\ga,-1)$   and consider the  operator $\bt{}= \ga b_1 -
>b_2$.
Observe that
\[
\ga s_k - t_k < 0 \iff \ga < \frac{t_k}{s_k}.
\]
Thus, if we write
\[
G_\ga = \{ k : \alpha < t_k/s_k\},
\]
then $p_{\tv,0}^-$ is the projection onto the subspace of $H$ spanned by
the corresponding
standard basis vectors and we have
\[
\tau(p_{\tv,0}^-) = \sum_{k \in G_\ga} \frac{1}{2^k}.
\]
Hence, if $\ga \ne \ga^\prime$ and we write $\tv^\prime = (\ga^\prime,0)$,
then $G_\ga
\ne G_{\ga^\prime}$, $\tau(p_{\tv,0}^-) \ne \tau(p_{\tv^\prime,0}^-)$ and so
\[
p_{\tv,0}^- \ne p_{\tv^\prime,0}^-.
\]

Hence there are an uncountable number of spectral projections of this form
and therefore
by part $(1)$ of Theorem 0.2, the set of extreme points of $B$ is uncountable.
\endexample

\bigskip
As the next two examples show, it is not always ``obvious" from the
geometry of  B whether
or not  $N$ is abelian.  Finding a condition on $B$ that is necessary and
sufficient for
commutativity of  $N$ is one of the major unsolved problems in the theory.

The pictures of the spectral scales in these examples are rather complicated.
The details are much clearer in color.  To view the pictures, go to

\centerline{http://www.math.psu.edu/anderson/pictures.}

\example Suppose that $H = L^2(0,1)$, $b_1$ is multiplication by
$x$, $b_2$ is the projection determined by the characteristic function of
the interval
$[0,1/2]$ and the trace  $\tau$ is integration. Thus, in this case $N$ is
abelian.
The spectral scale for these operators resembles two  bed sheets, tied at their
corners and billowing outwards from one another.  The coinciding edges of
the bed sheets
form the borders of 4 two dimensional faces  with parabolic boundaries that
meet at the
corners of the bed sheets.    These faces have the same shape as the
spectral scale for
multiplication by $x$.  (See \cite[end of \S1]{AAW}).  Thus, the spectral
scale of $b_1$
is visible in $B$.  The spectral scale of $b_2$ (which is just a
parallelogram) is also
visible  as the convex hull of the corners of the bed sheets.
\endexample

\bigskip

\example Now let us consider a non--commutative example.  Put $H
= L^2(0,1/2)\oplus L^2(1/2,1)$ and write
$$
b_1 =
\bmatrix
a&0\\0&a+1/2
\endbmatrix
\text{ and }
b_2 =
\bmatrix
1/2&1/2\\
1/2&1/2
\endbmatrix,
$$
where $a$ is multiplication by $x$ on $[0,1/2]$.   Thus, $b_1$ is still
multiplication by
$x$ on $[0,1]$, but the projection $b_2$ now does not commute with $b_1$.
In this case
the spectral scale is quite similar to the scale in the previous example,
although the
``billowing'' of the bed sheets is more pronounced.   So, based on this
evidence, there
does not seem to be an obvious visual method for determining when $N$ is
abelian in the
infinite dimensional case.
\endexample

\bigskip

 \section{Preliminary  Results on Convexity}
\medskip

How should one describe a face $F$ of a convex set $C$ in \rn{n}? Assuming
that we have
already some description of  $C$, it is natural to add some conditions to
that description
in order to distinguish the points of $F$ from the other points of $C$. The
most natural
conditions, from the convexity viewpoint, should be expressed in terms of
linear
functionals.  E.g. $F = \{x \in  C : f(x) = \alpha\}$. This only works if
the linear
functional  $f$  exposes the face  $F$.  Many faces are not exposed, so how
can we
describe them in some systematic way using linear functionals?  The answer
to that
question is to use more than one linear functional via the recursive
process described in
this section.

We  begin by reviewing some of the standard notions in the theory of
convexity and
then present some general results that will be used in the sequel. Although
some of the
material that appears here does not seem to be in any of the  standard texts on
convexity, it is quite possible that some of the unattributed results are
known.

\bigskip
 If $\tv = (t_1,\dots,t_n)$ is a nonzero vector in \rn{n} and $\ga$ is a
real number then
we write
\[
P(\tv,\ga) = \{\xv = ( x_1,\dots,x_n) \in \rnn{n}: t_1x_1 +\cdots + t_nx_n =
\lan\tv,\xv\ran = \ga\}
\]
for the {\bf hyperplane } in \rn{n} determined by $\tv$ and $\ga$.  The
half spaces of
\rn{n} determined by the hyperplane $P(\tv,\ga)$ are given as follows:
\[
P^\up(\tv,\ga) = \{\xv\in \rnn{n}: \lan\tv,\xv\ran \ge \ga\}, \text{ and }
P^\down(\tv,\ga) =
\{\xv\in \rnn{n}: \lan\tv,\xv\ran \le \ga\}
\]
Thus, $P(\tv,\ga)$ is the unique hyperplane in \rn{n} passing through $\xv$
and orthogonal
to
$\tv$ and so it is natural to call the vector $\tv$  a {\bf normal vector}
for this
hyperplane.

Now fix a  convex subset $C$ of \rn{n}. The hyperplane $P(\tv,\ga)$ is a
{\bf hyperplane
of support} for $C$ if there is a vector $\xv$ in the boundary of $C$ with
$\xv\in
P(\tv,\ga)$ and we have 
\[
C\subset P^\up(\tv,\ga).
\]
Note that in order for $P(\tv,\ga)$  to be a hyperplane of support we just
as easily
could  have required $C \subset P^\down(\tv,\ga)$.  In fact this is often
the definition
used in many texts on convexity.  We prefer to use the positive half space
since it is
more natural in our applications.

Note that  $P(\tv,\ga)$ is a hyperplane of support for $C$ containing the
vector $\xv \in
C$ if and only if
\[
\lan \xv'-\xv,\tv\ran \ge 0, \quad \text{for all }\xv' \in C.
\]

If $C$ is compact and $\tv$ is any nonzero vector in \rn{n}, then  it is
straightforward
to show that there is a unique real number $\ga$ such that $P(\tv,\ga)$ is
a hyperplane of
support for $C$.  Further, if $P(\tv,\ga)$ is a hyperplane of support for
$C$, then
\[
F(\tv,\ga) = P(\tv,\ga)\cap C
\]
is a face in the boundary of $C$.

 Let us now record a few general results on faces in convex sets. As above,
the symbol
$C$ stands for a convex set in \rn{n}.  The proof of the following Theorem is
essentially the same as the proof of Theorem 2.6.17 in \cite{Web}, although
its statement
is a bit more detailed than the one given there.

\begin{theorem}If  $P(\tv_1,\ga_1),\dots,P(\tv_k,\ga_k)$ are hyperplanes of
support for
$C$ such that
\[
C\cap P(\tv_1,\ga_1)\cap \cdots\cap P(\tv_k,\ga_k) = F(\tv_1,\ga_1)\cap
\cdots\cap
F(\tv_k,\ga_k)\ne \emptyset,
\]
$\gl_1,\dots,\gl_k$ are positive real numbers and we write

\begin{align*}
\ga_\gl &= \gl_1\ga_1 + \cdots \gl_k\ga_k \text{ and }\\
\tv_\gl &= \gl_1\tv_1+\cdots\gl_k\tv_k,
\end{align*}
then  $P(\tv_\gl,\ga_\gl)$ is a hyperplane of support for $C$ and
\[
F(\tv_\gl,\ga_\gl) =P(\tv_\gl,\ga_\gl)\cap C = F(\tv_1,\ga_1)\cap \cdots
\cap F(\tv_k,\ga_k).
\]
In particular, $F(\tv_\gl,\ga_\gl)$ is an exposed face of $C$.
\end{theorem}
\begin{proof} Fix $\xv \in F(\tv_1,\ga_1)\cap \cdots \cap F(\tv_k,\ga_k)$.
We have then
that
\[
 \lan\tv_i,\xv\ran = \ga_i
\]
for each $i$ and so
\[
\lan\tv_\gl,\xv\ran = \gl_1\lan\tv_1,\xv\ran + \cdots +
\gl_k\lan\tv_k,\xv\ran = \ga_\gl.
\]
 Hence $\xv \in F(\tv_\gl,\ga_\gl)$ and therefore,
\[
 F(\tv_1,\ga_1)\cap \cdots \cap F(\tv_k,\ga_k)\subset F(\tv_\gl,\ga_\gl).
\]

Now fix $\xv \in F(\tv_\gl,\ga_\gl)$ so that
\[
\lan\xv,\tv_\gl\ran = \ga_\gl = \gl_1\ga_1 +\cdots \gl_k\ga_k.
\]
If we had $\xv \notin F(\tv_i,\ga_i)$ for some $i$, then since $\xv \in C$
we would have
\[
\lan\tv_i,\xv\ran > \ga_i
\]
and therefore
\[
\lan\xv,\tv_\gl\ran = \gl_1\lan\tv_1\xv\ran + \cdots + \gl
_k\lan\tv_k,\xv\ran > \ga_\gl,
\]
contradicting $(*)$.  Hence,
\[ F(\tv_\gl,\ga_\gl) \subset F(\tv_1,\ga_1)\cap \cdots \cap F(\tv_k,\ga_k)
\]
and so
\[
F(\tv_\gl,\ga_\gl) = F(\tv_1,\ga_1)\cap \cdots \cap F(\tv_k,\ga_k).
\]

 Finally, since  $ F(\tv_1,\ga_1)\cap \cdots \cap F(\tv_k,\ga_k) \ne
\emptyset$ by
hypothesis, we get that $F(\tv_\gl,\ga_\gl)$ contains a point on the
boundary of $C$.
Further, if $\xv\in C$, then we have
\[
\lan\tv_i,\xv\ran \ge \ga_i\qquad i = 1,\dots, k
\]
and so $\lan\tv_\gl,\xv\ran \ge \ga_\gl$. Hence,
\[
C\subset P^\up(\tv_\gl,\ga_\gl),
\]
$P(\tv_\gl,\ga_\gl)$ is a hyperplane of support for $C$, and
$F(\tv_\gl,\ga_\gl)$
is  an exposed face of $C$.

\end{proof}
\medskip
\begin{remarks}

\noindent $(1)$ Observe that if $\mu > 0$, then $P(\tv,\ga) =
P(\mu\tv,\mu\ga)$. Thus,
If we let
\[
\mu = \frac{1}{\gl_1+ \cdots + \gl_k},
\]
and replace $\tv$ and $\ga$ with $\mu\tv$ and $\mu\ga$, then we have that
$\tv_\gl$ is a
convex combination of the $\tv_i$'s. Hence, the proposition is equivalent
to one that only
considers convex combinations of vector in \rn{n}.
\smallskip

\noindent $(2)$  If $F$ is a face in the boundary of $C$, then there is a
hyperplane of
support $P$ for $C$ that contains $F$ and such that $P$ is disjoint from
the relative
interior of $C$.  (See \cite[Theorem 11.3]{Rock}, for example).
\end{remarks}
\medskip
The following result is merely a restatement of Corollary 2.6.9 in
\cite{Web}

\begin{theorem} The intersection of any family of faces of $C$ can be
expressed as an
intersection of $n+1$ or fewer members of the family.
\end{theorem}

\begin{corollary} If $F$ is any nonempty face of $C$, then there is a (unique)
minimal exposed face $F_e$ such that
\[
F \subset F_e.
\]
\end{corollary}
\begin{proof} Let $\script{F}(F)$ denote the family of all exposed faces
that contain $F$.
We have that $\script{F}(F) \ne \emptyset$ by remark $(2)$ above.  Write
\[
F_e = \bigcap_{G\in\script{F}(F)} G.
\]
We have that $F_e$ is the intersection of a finite number of faces in
$\script{F}(F)$ by
Theorem 2.3.  Since $F \subset F_e$ this intersection is not empty and so
$F_e$ is an
exposed face by Proposition 2.1.  It is clear from its definition that
$F_e$ is the unique
minimal exposed face containing $F$.
\end{proof}

\begin{theorem} If $F$ is a face in a convex subset $C$ of \rn{n}, then
there is a unique
sequence
\[
F_1 \supset  \cdots \supset F_k
\]
Such that
\begin{enumerate}
\item $F_1$ is the unique minimal  exposed face in $C$ containing $F$.
\item If $1 <  i \le k$, then $F_i$ is the unique minimal exposed face of
$F_{i-1}$ that
contains $F$.
\item $F = F_k$.
\end{enumerate}
\end{theorem}
\begin{proof} The assertions are easily established using  Corollary 2.3
and a standard
induction argument.
\end{proof}

Let us call the chain $F_1 \supset  \cdots \supset F_k$ described above the
{\bf
minimal exposed facial chain} determined by $F$. 

\bigskip

\section{The Facial Structure of the Spectral Scale}
\medskip
With these preparations, we may now describe how  general (i.e.,
non--exposed) faces
of $B$ are determined by  spectral properties of certain operators in $N$.
It turns out
that these faces are described exactly as exposed faces are except that it
is necessary
to work in a ``cut--down" von Neumann algebra.

The key result that makes this analysis possible is the fact that a face
(of positive
dimension) in $B= B(b_1,\dots,b_n)$ may be identified as a spectral scale
in its own
right.  This scale is determined by   cut--down operators of the form
$rb_1r,\dots,rb_nr$ in the cut--down von Neumann algebra $rMr$.  In order
to describe the
situation more precisely, some new notation is required.

Suppose that $F$ is a  proper face of $B$ with dimension greater than zero;
i.e. $F$ is
not an extreme point.   In the next Lemma, it will be shown that $F$ is an
affine
translate of a certain related spectral scale. The new spectral scale is
determined  as
follows. By \cite[Theorems 1.1 and 2.4]{AAW} ,  we have $F
=\Psi([q^-,q^+])$ where  $q^-
\le q^+$ are projections in $N$.  Write $r =q^+-q^-$ and note that since
$F$ is a proper
face that is not an extreme point, $0  < r < 1$ and therefore $0 < \tau(r)
< 1$.  Next put
$M_r = rMr$, write
\[
\tau_r = \frac{1}{\tau(r)}\tau\big |M_r
\]
and define $\Psi_r$ on $M_r$ by the formula
\[
\Psi_r(x) = ( \tau_r(x), \tau_r(rb_1rx),\dots,  \tau_r(rb_nrx) ).
\]
Finally, let $B_F$ denote the spectral scale determined by $\Psi_r$
relative to $\tau_r$.

\begin{lemma} Suppose $F = \Psi([q^-,q^+])$ is a proper face in $B$ that is
not an extreme
point. Using the notation developed in the preceding paragraph, the map
$A$ defined by the formula
\[
A(\xv) = \frac{1}{\tau(r)}(\xv - \Psi(q^-))
\]
implements an affine isomorphism  of $F$ onto $B_F$.  Hence we have
\[
F = \Psi(q^-) + \tau(r)B_F = \Psi(q^-) + \Psi((M_r)_1^+) .
\]
\end{lemma}
\begin{proof}   Observe that we have the relation
\[
\tau(r)\Psi_r(x) = \Psi(x),\qquad x \in M_r.
\]

Now fix $\xv \in F$ so that we have $\xv = \Psi(a)$ for some $a \in
[q^-,q^+]$.  If we
write
\[
a_r = rar \text{ so that } a = q^- + a_r \text{ and } a_r \in (M_r)^+_1,
\]
then we have
\[
A(\xv) = A(\Psi(q^- + a_r)) =   \frac{1}{\tau(r)}(\Psi(q^-) +
\Psi(a_r)-\Psi(q^-)) =
\Psi_r(a_r). \phantom{0000} (*)
\]
and therefore $A(\xv) \in B_F$.  Further, since every vector in $B_F$ has
the form
$\Psi_r(a_r)$ for some $a_r \in (M_r)^+_1$, this same calculation shows
that $A$ is
surjective.

Finally, suppose $\xv_1$ and $\xv_2$ are vectors in $F$ such that $A(\xv_1)
= A(\xv_2)$.
As $\xv_1$ and $\xv_2$ are in $F$, we may find operators $a_1$ and $a_2$ in
$[p,q]$ such
that $\Psi(a_i) = \xv_i$ and we may write
\[
a_i = q^- + a_{i,r}, \quad i = 1,2,
\]
where  $a_{i,r} = ra_ir$.  With this we have
\[
A(\xv_1) =  \Psi_r(a_{1,r}) = \Psi_r(a_{2,r}) = A(\xv_2)
\]
by $(*)$. As
\[
\xv_i = \Psi(a_i) = \Psi(q^- + a_{i,r}),
\]
we get $\xv_1 = \xv_2$ and so $A$ is injective.  Hence, $A$ implements an
affine
isomorphism of $F$ onto $B_F$ and therefore
\[
F =  A^{-1}(B_F) = \Psi(q^-) + \tau(r)B_F.
\]
 \end{proof}
\bigskip


The next Lemma shows that if  B  is contained in one or more hyperplanes,
then  this
geometric condition is equivalent to the existence of linear dependence
relationships
(corresponding to the hyperplanes that contain $B$) among ${b_1, ..., b_n, 1}$.

\begin{lemma} If $0 \le k < n$, and the span of $B = B(b_1,\dots,b_n)$ has
dimension $k+1$, then  there are real numbers $s_1,\dots,s_{n-k}$ and linearly
independent vectors $\tv_1,\dots,\tv_{n-k}$ such that $\bt{i} = s_i1$ for $i =
1,\dots,n-k$.

Conversely, if there are real numbers $s_1,\dots,s_{n-k}$ and linearly
independent vectors $\tv_1,\dots,\tv_{n-k}$ such that $\bt{i} = s_i1$ for $i =
1,\dots,n-k$, then the subspace spanned by  $B$ has dimension at most $k+1$.
\end{lemma}
\begin{proof}  Suppose $P(\tv,s,\ga)$ is a hyperplane in \rn{n+1} that contains
$B$ and observe that since $0 \in B$, we may take $\ga = 0$.  In this case,
we have
\[
-s\tau(a)+\sum_{i=1}^n t_i\tau(b_ia) = \tau((\bt{}-s1)a) = 0
\]
for each $a\in M_1^+$.  Now select $\ga > 0$ such that $0 \le \ga(\bt{}-s1)^+
\le 1$.   Setting $a = \ga(\bt{}-s1)^+$ we get $\tau((\ga(\bt{}-s1)^+)^2) =
0$ and since
$\tau$ is faithful $(\bt{}-s1)^+ = 0$.  Similarly, we get $(\bt{}-s1)^- =
0$ and so
$\bt{}=s1$.

If $B$ has dimension $k+1$, then there are $n-k$ hyperplanes of the form
$P(\tv_i,s_i,0)$
that contain $B$ and such that their normal vectors $(-s_i,\tv_i)$ are linearly
independent.  In fact we must have that the $\tv_i$'s are linearly
independent.  To see
this, suppose that there are real numbers $\ga_2,\dots,\ga_{n-k}$ such that
\[
\tv_1 = \sum_{i=2}^{n-k}\ga_i\tv_i.
\]
In this case we get
\[
b_{\tv_1} = s_11 =  \sum_{i=2}^{n-k}\ga_i\bt{i} = \sum_{i=2}^{n-k}\ga_is_i 1.
\]
But, then $(-s_1,\tv_1)$ would be a linear combination of the remaining
normal vectors
contradicting their linear independence.

For the converse, if real scalars $\{s_1, ..., s_{n-k}\}$ and linearly
independent
vectors $\{\tv_1,\dots, \tv_{n-k}\}$  exist such that $\bt{i} = s_i1$  for
all $i = 1,
..., n-k$, then it immediately follows that every  $\xv$ in  $B$  is orthogonal
to each vector $(-s_i, \tv_i)$ and so the subspace spanned by $B$ has
dimension at most
$k+1$.
\end{proof}

\begin{corollary} Using the notation developed for Lemma 3.1 and letting $N_F$ be the von Neumann
algebra generated by $\{rb_1r, ..., rb_nr, r\}$, the following statements hold.
\begin{enumerate}

\item  $F$ has dimension one if and only if  $N_F$ has dimension one.

\item  If $F$ has dimension two, then $N_F$ is abelian.

\end{enumerate}
\begin{proof} We have $F = \Psi(q^-) + \Psi((M_r)_1^+)$.  Suppose that $N_F$ has dimension one,
i.e. that there are scalars $\{s_1, ..., s_n\}$ such that
$rb_ir = s_ir$ for each $i$.  In this case, if $\Psi(a) \in F$ so that  $a = q^- + a_r$, where $0
\le a_r \le r$, then we have
\begin{align*}
\Psi(a) &= (\tau(q^- + a_r), \tau(b_1(q^-+a_r),\dots,\tau(b_n(q^-+a_r))\\
& = \Psi(q^-) + (\tau(a_r), \tau(s_1a_r),\dots,\tau(s_na_r)\\
& = \Psi(q^-) + \tau(a_r)(1,s_1,\dots,s_n).
\end{align*}
Since $0 \le \tau(a_r)\le 1$, we get that $F$ is the line segment with end points
$\Psi(q^-)$ and $\Psi(q^-)+\Psi(r)$ and so $F$ has dimension one.

 If $F$ has dimension one, then by Lemma 3.1 $B_F$ has dimension one.  Thus by Lemma 3.2
there are $n$ linearly independent vectors
$\tv_1,\dots,\tv_n$ and real numbers $s_1,\dots,s_n$ such that $r\bt{i}r = s_ir$.  Since
the $\tv_i$'s are linearly independent, we get that each $rb_ir$ is a multiple of $r$. 
Hence assertion $(1)$ holds.

Now suppose $F$ has dimension two so that $B_F$ has dimension two by Lemma
3.1.  In this
case, by Lemma 3.2, there are $n-1$ linearly independent vectors $\tv_i$
and real
numbers $s_i$ such that $r\bt{i}r = s_ir$ for $i = 1,\dots, n-1$.  This means
that
there is a vector $\tv$  such that
\[
\spn(r\bt{}r,1) = \spn(rb_1r,\dots,rb_nr,1)
\]
and so $N_F$ is abelian.
\end{proof}
\end{corollary}

\remark

Observe that the spectral scale described in Example 1.3 has two
dimensional faces
such that $M_r$ is abelian (in fact $M$ is abelian in this case), but $M_r$
is infinite
dimensional since it is isomorphic to the von Neumann algebra generated by
multiplication
by $x$ on a suitable interval.

\bigskip

Let us return now to the  description of all of the faces of $B$. By
identifying a
face $F$ with $B_F$, we may view all subfaces of $F$ as faces in $B_F$. In
particular, we
may apply the results in Theorem 0.2 to the exposed subfaces of $B_F$ to
obtain data
about the corresponding subfaces of $F$.  This process involves a good deal
of new
notation and is a bit complex.  We now gather this information in a definition.

\definition In what follows, we shall consider cut--downs of  operators
in $M$ of the form
\[
b^\prime = rbr,
\]
where $r$ is a projection and $b$ is self-adjoint.  Note that we may view
$b^\prime$
as an element of $M$ or the von Neumann algebra $rMr$. If  $s$ is a real
number then the
{\bf cut down spectral interval projections} determined by $b, r$ and $s$
are the
spectral projections of $b^\prime$  determined by the intervals $(-\infty,s)$
and  $(-\infty,s]$.  We compute the spectrum and the corresponding spectral
projections
in the algebra $M_r$ where the projection $r$ is the identity.

If $\spk{1},\spk{2},\dots,\spk{k}$ is a sequence of spectral pairs, then
its associated
{\bf facial complex} is defined inductively as follows.  Write
\[
r_1 = \ptspp{1}-\ptsmm{1} \text{ and } \bt{2}^\prime = r_1\bt{2}r_1.
\]
Let $\qtsppmm{2}$ denote the associated cut down spectral interval
projections determined
by $\bt{2},r_1$ and $s_2$ and write
\[
r_2 = \qtspp{2}-\qtsmm{2}.
\]

Now suppose that for some $1 < i < k$ we have defined sequences

\begin{align*}
&r_1 > r_2 > \cdots > r_i,\\
&\bt{2}^\prime,\dots, \bt{i}^\prime,\text{ and }\\
&\qtsppmm{2},\dots,\qtsppmm{i}
\end{align*}
such that if $2 \le j\le i$, then
\[
 \bt{j}^\prime  = r_{j-1}\bt{j}r_{j-1}, \quad r_j = \qtspp{j}-\qtsmm{j}
\]
where $\qtsppmm{j}$ are the cut down spectral interval projections
determined by
$\bt{j},r_{j-1}$ and $s_j$.

We may then set
\[
\bt{i+1}^\prime = r_i\bt{i+1}r_i,
\]
let $\qtsppmm{i+1}$ denote the spectral interval projections determined by
$\bt{i+1},r_i$
and $s_{i+1}$ and write
\[
r_{i+1} = \qtspp{i+1}-\qtsmm{i+1}.
\]
This completes the induction.

Thus the facial complex determined by the spectral pairs
$\spk{1},\dots,\spk{k}$ consists
of the projections $\ptsppmm{1}$ together with the sequences

\begin{align*}
&r_1 > \cdots > r_{k-1},\\
&\bt{2}^\prime,\dots,\bt{k}^\prime,\text{ and }\\
&\qtsppmm{2},\dots,\qtsppmm{k}
\end{align*}

as defined above.
\enddefinition

\begin{theorem} If $F$ is a face in the spectral scale and $F_1 \supset \cdots
\supset F_k = F$ is its associated minimal exposed facial chain, then there
is a sequence
$\spk{1},\spk{2},\dots,\spk{k}$ of spectral pairs such that if
$\qtsppmm{2},\dots,\qtsppmm{k}$ denote the spectral interval projections in its
associated facial complex, then $F_1 = \Psi([\ptsmm{1},\ptspp{1}])$ and for
$2 \le i \le
k$ we have
\[
F_i =  \Psi([q_i^-,q_i^+]) = \Psi(\qtsmm{i-1}) +\Psi([\qtsmm{i} ,\qtspp{i}])
\]
where
\[
q_i^- = \ptsmm{1} + \qtsmm{2} + \cdots + \qtsmm{i}\text{ and } q_i^+ =
\ptsmm{1} +
\qtsmm{2} + \cdots + \qtsmm{i-1}+\qtspp{i}.
\]
\end{theorem}
\begin{proof}  Let us proceed by induction.  As $F_1$ is an exposed face in
$B$, there is
a spectral pair $\spk{1}$ such that
\[
F_1 = \Psi([\ptsmm{1},\ptspp{1}]).
\]
If $F = F_1$, we are done.

Otherwise,  write $r_1 = \ptspp{1}-\ptsmm{1}$ and let $A_1$ denote the
affine map defined
in Lemma 3.1 for the face $F_1$ of $B$. We have that $A_1(F_2)$ is an
exposed face in
$B_{F_1}$ and so  by Theorem 0.2, there is a spectral pair$\spk{2}$ such
that if
$\bt{2}^\prime = r_1\bt{2}r_1$ and
$\qtsppmm{2}$ denote the associated cut down spectral interval projections,
then we have
\[
A_1(F_2) = \Psi_{r_1}([\qtsmm{2},\qtspp{2}]).
\]
Hence, by Lemma 3.1 we have
\[
F_2 = A_1^{-1}(\Psi_{r_1}([\qtsmm{2},\qtspp{2}])) = \Psi(\ptsmm{1}) +
\Psi((M_{r_1})_1^+)
\]

If $k = 2$, then we are done.  Otherwise, suppose that spectral pairs
$\spk{1},\spk{2},\dots,\spk{i}$ and sequences

\begin{align*}
&r_1 > \dots > r_{i-1},\\
&\bt{2}^\prime,\dots, \bt{i}^\prime,\text{ and }\\
&\qtsppmm{2},\dots,\qtsppmm{i}
\end{align*}

have been determined as above such that if $2 \le j < i$ and
\[
q_j^- = \ptsmm{1} + \qtsmm{2} + \cdots + \qtsmm{j}\text{ and } q_i^+ =
\ptsmm{1} +
\qtsmm{2} + \cdots + \qtsmm{j-1}+\qtspp{j},
\]
then we have
\[
F_j = \Psi([q_j^-,q_j^+]) = \Psi(\qtsmm{j-1}) + \Psi([\qtsmm{j},\qtspp{j}]).
\]

With this, we may write $r_{i} = \qtspp{i}-\qtsmm{i}$, $\bt{i+1}^\prime =
r_i\bt{i+1}r_i$
and let $\qtsppmm{i+1}$ denote the associated spectral interval
projections.   As $F_{i}$ is exposed in $F_{i-1}$, we may apply Lemma 3.1
to get the
affine map $A_{i-1}$ such that $A_{i-1}(F_{i-1}) = B_{F_{i-1}}$ and
$A_i(F_{i})$ is an
exposed face in $B_{F_{i-1}}$. Hence, applying Theorem 0.2 we get a
spectral pair
$\spk{i}$ such that if we write $\bt{i}^\prime = r_{-1}\bt{i}r_{i-1}$ and let
$\qtsppmm{i}$ denote the associated spectral interval projections, then we have
\[
A_{i-1}(F_{i}) = \Psi_{r_{i-1}}([\qtsmm{i},\qtspp{i}]).
\]
Applying $A_{i-1}^{-1}$ to $A_{i-1}(F_{i})$, yields the desired result.
\end{proof}
\bigskip

\section{Corners and their Generalizations}
\medskip
In the case of a single operator the corners on the boundary of the
spectral scale are in
one-to-one correspondence with the gaps in the spectrum of $b$. In this
section we study
analogous faces in higher dimensions.

\definition If $C$ is a convex set in a vector space and $F$ is a face in
$C$, then we say $F$ is a {\bf sharp} face if it is contained in two or
more hyperplanes
of support.
\enddefinition

Recall that  if $A$ is an affine space in \rn{n}, then  the dimension of
$A$ is the
dimension of $A-\xv$, where $\xv$ is any vector in $A$.  The  {\bf
dimension} of the
convex set $C$ is by definition the dimension of the affine space that $C$
generates.
The dimension is denoted by $\dim(C)$.

Now fix a nonempty face $F$ of $C$. Recall that a vector $\tv\in \rnn{n}$
is a {\bf normal
vector} for $F$ if $P(\tv,\ga(\tv))$ is a hyperplane of support for $C$ and
$C\subset
P^\up(\tv,\ga(\tv))$. The {\bf normal cone} of $F$ is by definition
\[
K_F = \{\tv\in\rnn{n}: \tv\text{ is a normal vector for $F$}\}.
\]

The following result is  basic to our analysis of sharp faces in the
spectral scale.
\smallskip
\begin{theorem} If $F = \Psi([q^-,q^+])$ is a face in the spectral scale
$B$ and
\spe is a spectral pair,  then the following statements are equivalent.

\begin{enumerate}
\item  The vector $(-s,\tv)$ is in $K_F$.
\item We have  $\ptsm\le q^- \le q^+\le \ptsp$.
\item If we write
\[
r_1 = q^-, r_2 = q^+-q^- \text{ and } r_3 = 1-q^+,
\]
set $ b_{\tv,1} = \bt{}r_1$ and $ b_{\tv,2}=\bt{}r_3 $ , then

\begin{align*}
b_{\tv,1} &\le sr_1,\\
 \bt{}r_2 &= sr_2,\\
 b_{\tv,2}  &\ge sr_3,
\end{align*}

 and $q^\pm\bt{} =\bt{}q^\pm$.

In other words,  the matrix for $\bt{}$ determined by  the projections
$r_1,r_2$ and
$r_3$ has the form
\[
\bt{} =
\bmatrix
b_{\tv,1}&0&0\\
0&s&0\\
0&0&b_{\tv,2}
\endbmatrix.
\]
In particular, the $r_i$'s commute with $\bt{}$.
\end{enumerate}
\end{theorem}

\begin{proof} Suppose $(1)$ holds and $P(s,\tv,\ga)$ is the hyperplane of
support that
contains $F$. Since $F = \Psi([q^-,q^+]) \subset \Psi([\ptsm,\ptsp])$  we get
\[
\ptsm\le q^- \le q^+\le \ptsp
\]
by parts $(2)$ and $(3)$ of Theorem 0.2.  Hence $(1)\implies (2)$. It is
also clear from
the above that $(2)\implies (1)$.

Now suppose $(2)$ holds.  In this case since, $\ptspm$ are the spectral
interval
projections for $\bt{}$ and $\ptsm\le q^- \le q^+\le \ptsp$ we get
\[
\bt{}q^- \le sq^-, (q^+-q^-)\bt{} = s(q^+-q^-) \text{ and } \bt{}(1-q^+)
\ge s(1-q^+).
\]
Further we have $q^-\bt{} = \bt{}q^-$ and $(1-q^+)\bt{} = \bt{}(1-q^+)$ and
so $\bt{}$ has
the desired form.  Thus, $(2)\implies (3)$.

Now suppose $(3)$ holds.  As $b_{\tv,1} \le sr_1, \bt{}r_2 = sr_2 \text{
and } b_{\tv,2}
\ge sr_3$, we get
\[
\ptsm \le q^- < q^+ < \ptsp
\]
and so $(3)\implies (2)$.

\end{proof}
\medskip
\begin{theorem} If $F = \Psi[q^-,q^+])$ is a face in the spectral scale,
then $F$
is a sharp face if and only if there are linearly independent spectral
pairs $(s_1,\tv_1)$
and $(s_2,\tv_2)$ such that
\[
F \subset P(s_1,\tv_1,\ga_1)\cap P(s_2,\tv_2,\ga_2).
\]
In this case if we write
\[
r_1 = q^-,\quad r_2 = q^+-q^-,\text{ and } r_3 = 1-q^+,
\]
then we have
\begin{enumerate}
\item $r_1+r_2+r_3 = 1$,
\item $\bt{i}r_1 \le s_ir_1$,
\item $\bt{i}r_2 = s_ir_2$, and
\item $\bt{i}r_3 \ge s_ir_3$
\end{enumerate}
for $i = 1,2$ and the $r_i$'s commute with $\bt{}$.
\end{theorem}
\begin{proof} The assertions are easily established using the definition of
a sharp face
and Theorem 4.2.
\end{proof}
\medskip
\begin{theorem} If $F = \Psi([q^-,q^+])$ is a proper face of $B$ and
$F_e$ is the minimal exposed face containing $F$, so that $F_e =
\Psi([p_F^-,p_F^+])$, where $p_F^\pm$ are spectral interval  projections of
a linear combination of the $b_i$'s,
then  we have
\[
p_F^- \le q^-\le q^+ \le p_F^+
\]
Further, if \spe is a spectral pair such that   $\Psi([\ptsm,\ptsp])$ is an
exposed face
of $B$ that contains $F_e$, then
\[
\ptsm \le p_F^- \le q^- \le q^+ \le p_F^+ \le \ptsp.
\]
\end{theorem}
\begin{proof} All of the assertions are simple consequences of part $(1)$
and $(2)$ of
Theorem 4.2.
\end{proof}
\medskip
\begin{corollary} Using the notation of Theorem 4.4, a face $F =
\Psi([q^-,q^+])$  of $B$
is an extreme point of
$B$ if and only if $q^- = q^+$.  In this case we have
\[
\bt{} =
\bmatrix
b_{\tv,1}&0\\
0&b_{\tv,2}
\endbmatrix,
\]
where $b_{\tv,1} \le sq^- $ and $b_{\tv,2} \ge s(1-q^+)$.   The face $F =
\Psi(q^\pm)$ is
an exposed point if and only if $p_{F^-} = p_{F^+}$.
\end{corollary}
\begin{proof}   We have that $F$ is a singleton if and only if $q^- = q^+$.
In this case, in the notation of Theorem 4.2, we get $r_2 = 0$ and so
$\bt{}$ has the
desired matrix decomposition by part $(3)$ of  Theorem 4.2.  We have that
$\xv$ is an
exposed point if and only if  $F = F_e$, where $F =\{\xv\}$.  Since $F_e =
\Psi([p_F^-,p_F^+])$ by Corollary 4.4,  the final assertion follows.
\end{proof}

Before presenting the next result, it is convenient to introduce some new
notation.
Suppose $F$ is a face in a convex subset of \rn{n}.  The {\bf degree} of
$F$ is denoted
by $\deg(F)$ and is  by definition the  dimension of the normal cone $K_F$.
Recall
that, the subspace spanned by $K_F$  is $K_F - K_F$ and so the degree of
$F$ is just the
dimension of $K_F-K_F$.

\begin{theorem}The following assertions hold.
\begin{enumerate}
\item The degree of $F$ is equal to the maximal number of linearly
independent vectors in
$K_F$.
\item We have $\deg(F) + \dim(F) \le n$.
\item $F$ is a sharp face if and only if $\deg(F) > 1$.
\end{enumerate}
\end{theorem} \begin{proof} For $(1)$, it is clear that there can be no
more than
$\deg(F)$ such vectors.    Now suppose $\{\tv_1^+,\dots ,\tv_k^+\}$ is a
maximal set
of linearly independent vectors in $K_F$ and fix any vector $\tv$ in
$K_F-K_F$  We must
have $\tv = \tv^+ -\tv^-$, where $\tv^\pm \in K_F$.  This means that both
$\tv^+$ and
$\tv^-$ are linear combinations of $\{\tv_1^+,\dots ,\tv_k^+\}$ and so
$\tv$ is a
linear combination of $\{\tv_1^+,\dots, \tv_k^+\}$. Thus, $\{\tv_1^+,\dots,
\tv_k^+\}$ is
a basis for $K_F-K_F$.

Next, since each vector in $K_F$ is orthogonal to the difference of any two
vectors in
$F$, assertion $(2)$ is clear. Finally, assertion $(3)$ follows immediately
from the
definitions of a sharp face and the degree of $F$.

\end{proof}

\begin{corollary}If $B \subset \rnn{n+1}$,  $F= \Psi([q^-,q^+])$ is a face
in $B$  and there
are spectral pairs
$\spk{1},\dots, \spk{n}$  such that
\begin{enumerate}
\item $F \subset \Psi([\ptsmm{i},\ptspp{i}])$ for $i = 1,\dots, n$ and
\item $\tv_1,\dots,\tv_n$ are linearly independent,
\end{enumerate}
then  $q^-$ and $q^+$ are central projections in the von Neumann algebra
generated by
$\{b_1,\dots,b_n\}$.
\end{corollary}
\begin{proof} Applying  Theorem 4.3, we get that  $q^-$ and $q^+$ commute
with each
$b_{\tv_i}$.  Since $\tv_1,\dots,\tv_n$ are linearly independent, each
$b_i$ is a linear
combination of the $b_{\tv_i}$'s and so $q_1$ and $q^+$ commute with each
$b_i$.
\end{proof}

\remark
\noindent  In the notation of the
last Corollary, $\deg(F) \ge n$, so that
 by part $(2)$ of Theorem 4.6, we get that $\dim(F) \le 1$.  Thus, in
this case, $F$ is either a sharp extreme point or a sharp line segment in $B$.
\endremark

\begin{corollary} If $F = \Psi([q^-,q^+])$ is a face in the spectral scale and
$(-s_1,\tv)$ and $(-s_2,\tv)$ are distinct vectors in $K_F$ with $s_1 <
 s_2$, then
\begin{enumerate}
\item  The interval $(s_1,s_2)$ lies in a gap in the spectrum of $\bt{}$ and
\item $F$ is an exposed  point.
\end{enumerate}
\end{corollary}
\begin{proof} Applying Theorem 4.3 (and using its notation) we get that if
$r_2 =
q^+-q^-$, then
\[
\bt{}r_2 = s_1r_2 = s_2r_2
\]
and so $r_2 =  0$  and $q^-= q^+ = q$.

Next as $(-s_i,\tv) \in K_F$ for $i = 1,2$, we get
\[
\ptsmm{i} \le q \le \ptspp{i}
\]
for $i = 1,2$.  On the other hand, since $s_1 < s_2$ we have $\ptspp{1} \le
\ptsmm{2}$ so
that
\[
\ptsmm{2} \le q \le \ptspp{1} \le \ptsmm{2}
\]
and therefore
\[
q = \ptspp{1} = \ptsmm{2}.
\]
Hence there must be a gap in the spectrum of $\bt{}$ that contains the interval
$(s_1,s_2)$.

Next, if we select any $s$ such that $s_1 < s < s_2$, then we have $\ptsm =
\ptsp = q$ and
so $F = \Psi(\ptspm)$.  Thus, $F$ is an  exposed  point by Corollary 0.3.
\end{proof}

\section{Isolated  extreme points and central projections}
\medskip

We say that an element $\xv$ in a convex set $C$ is an {\bf isolated
extreme point} of
$C$ if it is an extreme point and it is isolated in the set $\e(C)$ of all
extreme points
of $C$. Our next goal is to show that isolated extreme
points in $B$ are images of central projections in $N$ and further that
this condition
characterizes isolated extreme points in the finite dimensional case.  It
is convenient
to begin by presenting three Lemmas. The first Lemma simply records an easy
fact about
convex subsets in \rn{n}.  The (slick) proof of the second is due to Robert
Phelps
(private communication).
\bigskip

Recall that if $A_1$ and $A_2$ are disjoint, compact, convex subsets in a
real topological
vector space, then there is a continuous linear functional $f$ that {\bf
separates} $A_1$ and $A_2$ in the sense that  there are real real numbers
$s < t$ such
that
\[
f(\xv_1) < s < t < f(\xv_2)
\]
for all $\xv_1\in A_1$ and $\xv_2\in A_2$ \cite[Theorem 3.4 (b)]{Rud:a} .
In this case we
say that $s$ and $t$ are {separation constants} for $f$.

\begin{lemma} If $C_0$ is a compact, convex subset of \rn{n}, $\xv$ is a
point in
\rn{n} disjoint from $C_0$ and $C = \conv(C_0,\xv)$, then a linear functional
$f$ exposes $\xv$ in $C$ if and only if $f$ separates $C_0$ and $\{\xv\}$.
\end{lemma}
\begin{proof} Suppose that $f$ exposes $\xv$ in $C$ so that there is a
number $t$ such
that $f(\xv) = t$ and $f(\zv) < t$ if $\zv \in C$ and $\zv \ne \xv$.  Write
\[
r = \sup\{f(\yv): \yv \in C_0\}
\]
and suppose $\yv_n$  is a sequence in $C_0$ such that $\lim_n f(\yv_n) =
r$.  As $C_0$ is
compact, we may assume that $\yv_n$ converges to a vector $\yv \in C_0$.
Since $f(\yv) =
r$, $f$ exposes $\xv$ and $\yv \ne \xv$, we must have that $r < t$.  Hence,
we get that
$f$ separates $C_0$ and $\xv$.

Now suppose that $f$ separates $C_0$ and $\xv$ with
separation constants $s < t$ and fix
$\zv \ne \xv \in C$.  Since  $C = \conv(C_0,\xv)$,  there is a number $r$
with $0 \le r <
1$ and a vector $\yv \in C_0$, such that  $\zv = r\xv + (1-r)\yv$.  Hence,
$f(\zv) < rt +
(1-r)s < t$ and so $f$ exposes $\xv$.

\end{proof}
\noindent \begin{lemma} If $C$ is a convex, compact subset in \rn{n}
and  $\xv$ is an isolated extreme point in $C$,  then the set of linear
functionals that
expose $\xv$ is a nonempty open set.  Further, there are linearly
independent vectors
$\tv_1,\dots,\tv_n$ and real numbers $\ga_1,\dots,\ga_n$ such that the
hyperplanes
$P(\tv_1,\ga_1),\dots,P(\tv_n,\ga_n)$ are hyperplanes of support for $C$ and
\[
P(\tv_i,\ga_i)\cap C = \{\xv\},\quad i = 1,\dots,n.
\]
\end{lemma}
\begin{proof} Write  $A$ for the closure of $\e(C)\setminus \{\xv\}$ and
observe that
since  $\xv$ is an isolated point in $\e(C)$, $\xv \notin A$.   Further, if
$C_0$
denotes the closed convex hull of $A$, then $\xv \notin C_0$.  Indeed, if
we had
$\xv\in C_0$, then $\xv$ would be an extreme point of $C_0$.  But in this
case we would
get that $\xv \in A$ by Choquet's ``converse" to the Krein--Milman theorem
\cite[Appendix
B14]{Dix}.

As noted in the paragraph before Lemma 5.1, there is a continuous linear
functional $f$
that separates $C_0$ and $\xv$ and so, by Lemma 5.1 $f$ exposes $\xv$ in
$C$.  Hence the
set of linear functionals that expose $\xv$ is not empty.

Now suppose that $f$ is a linear functional that exposes $\xv$.  Applying
Lemma 5.1, we
get that there are real numbers $s<t$ such that  $f$ separates $C_0$ and
$\xv$ with separation constants $s$ and $t$.   Now select $\eps > 0$ so that if
\[
M = \sup\{\|\yv\|:\yv \in C_0\}, \text{ then } \max\{\eps\|\xv\|, \eps M \} <
\frac{t-s}{3}
\]
and let $g$ denote a linear functional with $\|f-g\| < \eps$.  With this we get
\[
g(\xv) = f(\xv) + g(\xv) - f(\xv) >  t -\eps\|\xv\| \ge t - \frac{t-s}{3}.
\]
If $\yv \in C_0$, then we have
\[
g(\yv) = f(\yv) + g(\yv) - f(\yv) < s + \eps M \le s + \frac{t-s}{3}
\]
and so $f$ separates $C_0$ and $\xv$ with separation constants  $s +
(t-s)/3$ and
$t- (t-s)/3$.  Applying  Lemma 5.1 we get that $g$ also exposes $\xv$ and
so the set of
all such linear functionals is open.

To get the final conclusion, observe that since the set of linear
functionals that expose
$\xv$ is a nonempty open set, it contains a linearly independent set of $n$
linear
functionals,  $f_1, \dots, f_n$.  Thus there are linearly independent vectors
$\tv_1,\dots,\tv_n$ such that $f_i$ is given by the inner product with
$\tv_i$, and the
conclusion follows.
\end{proof}

\remark Observe that the proofs of the last two lemmas go through in any
normed linear
space.
\endremark

\noindent \begin{theorem} The following statements hold.
\begin{enumerate}
\item If $\xv$ is an isolated extreme point of $B$, then there is a central
projection
$p\in N$ such that $\xv = \Psi(p)$.

\item If $N$ is finite-dimensional and $\xv$ is an extreme point of $B$,
then  $\xv =
\Psi(p)$ for some central projection $p \in N$, if and only if  $\xv$ is an
isolated
extreme point of
$B$.
\end{enumerate}
\end{theorem}

\begin{proof} (1) Suppose $\xv$ is an isolated extreme point in $B$.  By
part $(1)$ of
Theorem 0.2, there is a unique projection $p\in N$ such that $\xv =
\Psi(p)$.  Next,
applying Lemma 5.2 and parts  $(2)$ and  $(3)$ of Theorem 0.2, we get
linearly independent
spectral pairs$\spk{1},\dots,\spk{n+1}$ such that
\[
\{\xv\} = \Psi([\ptsmm{k},\ptspp{k}]),\quad k = 1,\dots,n+1.
\]
Since $\Psi^{-1}(p) = \{\xv\}$, we have
\[
p  =\ptsppmm{k}\quad k = 1,\dots,n+1
\]
Since these spectral pairs are linearly independent, a straightforward
argument shows
that the subspace of \rn{n} spanned by $\tv_1,\dots,\tv_{n+1}$ has
dimension $n$.  Hence, relabeling if necessary, we may assume that
$\tv_1,\dots,\tv_n$ are
linearly independent. Therefore, by Corollary 4.7 $p = \ptsppmm{k}$ is a
central
projection in $N$.

(2)  Suppose that $\xv = \Psi(p)$, where $p$ is a central projection in $N$
and let
$\xv_m$ denote a sequence of extreme points of $B$ that converges to
$\xv$.  Applying part $(1)$ of Theorem 0.2  we get a sequence $p_m$ of
unique projections
in $N$ such that $\Psi(p_m) = \xv_m$ for each $m$.  Since $M$ is finite
dimensional, we
may pass to a  subsequence and assume that the sequence $\{p_m\}$ is norm
convergent with
limit $q$.  Since $\Psi$ is continuous, we have  $\Psi(q) = \xv = \Psi(p)$
and since
$\Psi^{-1}(\xv) = \{p\}$, we get $q = p$.  Now, since $p$ is central and
the central
projections are isolated in the set of all projections  in $N$, we must
have that
$p_m = p$  (and hence $\xv_m = \xv$) for
all sufficiently large  $m$.  Hence $\xv$ is an isolated extreme point.

If $\xv$ is an isolated extreme point in $B$, then by part $(1)$ of this
Theorem, there
is a central projection $p$ in $N$ such that $\xv = \Psi(p)$

\end{proof}
\bigskip

\remarks

\medskip
\noindent $(1)$ Let $b$ denote multiplication by $x$ on the interval
$[0,1]$ and write $N
= \{b\}^{\prime\prime} =  L^\infty(0,1)$. In this case every  projection in
$N$ is central
but $B$ has no isolated extreme points. (See the example at the end of
section 1 in
\cite{AAW}). Hence, the converse to part $(1)$ of Theorem 5.3 is false.

\medskip
\noindent $(2)$ If $b,B$ and $N$ are as in remark $(1)$ above, then there
are linearly
independent spectral pairs $\spk{1}$ and $\spk{2}$ such that
\[
\{0\} = \Psi([\ptsmm{i},\ptspp{i}]),\quad i = 1,2,
\]
but $0$ is not an isolated extreme point in $B$. (See \cite{AAW}).  Thus,
the converse to Lemma 5.2 is also false.
\endremark
\smallskip

\noindent \begin{corollary}  $N$ is abelian and has finite dimension if and
only if $B$
has a finite number of extreme points.
\end{corollary}

\begin{proof} If $N$ is abelian and has finite dimension, then  each
projection in $N$ is
diagonal.  Hence $N$ contains a finite number of projections and since
every extreme
point in $B$ is the image under $\Psi$ of a projection $N$, $\e(B)$ is finite.

For the converse, if  $B$ has a finite number of extreme points,   then
each extreme
point in $B$ is an isolated extreme point and so by part $(1)$ of Theorem
5.3, each
extreme point of $B$ has the form $\Psi(p)$, where $p$ is a central projection.
On the other hand, if $\spk{}$ is any spectral pair, then $\Psi(\ptspm)$ is
an extreme
point by part $(1)$ of Theorem 0.2 and so $\ptspm$ is a central projection.
Thus, for
each fixed $\tv$ the von Neumann algebra
\[
\{\ptspm: s\in \mathbb R\}^{\prime\prime}
\]
is central and therefore each  $\bt{}$ is central.  Hence $N$ is abelian.
Since the
projections of the form $\ptspm$ generate $N$, this algebra  has  finite
dimension.

\end{proof}
\bigskip

\section{Countable extreme points}
\medskip

Our goal in this section is to establish  the following Theorem.
\bigskip

\begin{theorem} If the set $\e(B)$ of extreme points of  $B$ is countable,
then $N$ is
abelian.
\end{theorem}

It is convenient to present the bulk of  the proof in a series of Lemmas.
We begin with
a simple topological fact.

\begin{lemma} If $X$ is a countable locally compact Hausdorff space, then the
isolated points in $X$ are dense.
\end{lemma}
\begin{proof}  Suppose that the closure of the isolated points of $X$ is
not all of $X$.
 By \cite[Theorem 2.7]{Rud:b}  each open subset of  $X$ is a locally
compact Hausdorff
space.   Hence, there is a nonempty open set $U$ in $X$ with no isolated
points.  Since
$X$ is countable, $U$ is countable and we may write $U = \{x_1,
x_2,\dots\}$.  Set $U_n =
U\setminus\{x_n\}$ for each $n$.  As $X$ is Hausdorff, each $U_n$ is open
and since $x_n$
is not isolated, each $U_n$ is dense in $U$. Hence by the Baire Category
Theorem
\cite[Theorem 2.2 (b)]{Rud:a} , $\cap U_n$ is dense in $U$.  Since, $\cap
U_n =
\emptyset$, this is impossible. Hence, no such $U$ exists and so the
isolated extreme
points in $X$ are dense.
\end{proof}

\begin{lemma} It suffices to prove the theorem in the case where $n = 2$.
\end{lemma}
\begin{proof} Write $N_{ij} = \{b_i,b_j,1\}''$, let $\pi_{ij}$  denote the
map from $B$ to
\rn{3} defined by the formula
\[
\pi_{ij}((\tau(a),\tau(b_1a),\dots,\tau(b_na))=
(\tau(a),\tau(b_ia),\tau(b_ja)), \quad
a\in N_1^+
\]
and write $B_{ij} = \pi_{ij}(B)$.  It is clear that $B_{ij} = B(b_i,b_j)$.

Next suppose $\e(B)$ is countable and fix $i$ and $j$.  For each extreme
point $\xv$ in
$B_{i,j}$  write $F_{\xv} = \pi_{ij}^{-1}(\xv)$ so that $F_\xv$ is a closed
face of $B$.
Now each $F_\xv$ is nonempty and so contains extreme points.  Since
$F_{\xv_1}$ and
$F_{\xv_2}$ are disjoint if $\xv_1$ and $\xv_2$ are distinct extreme
points, the extreme
points of $B_{ij}$ must be countable.
Thus, if $\e(B)$ is countable, then $\e(B_{ij})$ is countable for each pair
$(i,j)$.
Hence, if the Theorem is true when $n = 2$, then each $N_{ij}$ is abelian
and therefore
$N$ is abelian.
\end{proof}

For the remainder of this section, we shall assume that $N$ is the von
Neumann algebra
generated by $b_1,b_2$ and the identity and that $B = B(b_1,b_2)$ denotes
the associated
spectral scale. As $B \subset \rnn{3}$ each face in $B$ must have dimension
zero, one or
two. Recall that for each real $s$, with $0\le s \le 1$
\[
I_s = \{(y,z)\in  R^2: (s,y,z) \in B\}
\]
denotes the {\bf isotrace slice} of $B$ at $s$.
\begin{lemma} If $F$ is a face in $B$,
then the following statements hold.
\begin{enumerate}
\item If $F$ has dimension greater than zero, then $F$ is transverse to the
isotrace
slices of $B$.

\item If $F$  has dimension one, then $F$ meets each isotrace slice in at
most one point.

\item If $F$ has dimension two, and $I_s$ contains a point in the relative
interior of
$F$, then $F\cap I_s$ is a proper line segment.
\end{enumerate}
\end{lemma}
\begin{proof}
We have that $F = \Psi([\ptsm,\ptsp])$ for some  spectral pair $(s,\tv)$ By
part $(3)$ of
Theorem 0.2. Now suppose that $F$ has dimension greater then zero. In this
case we must
have  $\ptsm < \ptsp$ and so $\tau(\ptsm)< \tau(\ptsp)$.    Hence, $F$ is
transverse to
the isotrace slices of $B$ and so $(1)$ is true.

The remaining assertions now follow easily from $(1)$.
\end{proof}

\begin{lemma}If $0 \le s \le 1$ then the  following statements hold for $I_s$.
\begin{enumerate}
\item If $G$ is a proper closed line segment in the boundary of $I_s$, then
there is a two dimensional face $F$ of $B$ such that $G \subset F\cap I_s$.

\item  If $\xv$ is an extreme point of $I_s$, then either $\xv$ is an
extreme point of
$B$ or else it lies in the relative interior of a one dimensional face of $B$.

\item If $\xv$ is an extreme point in $I_s$ that is an isolated point in
$\e(I_s)$, then there are two distinct two dimensional faces of $B$ that
contain $\xv$.

\item If the set of extreme points $B$ is countable then the set of extreme
points of
$I_s$ is countable.
\end{enumerate}
\end{lemma}
\begin{proof} Fix a closed proper line segment $G$ in the boundary of
$I_s$. As $G$ is in
the boundary of $I_s$, it is in the boundary of $B$ and so is contained in
a minimal face $F$ in $B$.  Since faces in $B$ of dimension greater than
zero are
transverse to the isotrace slices by Lemma 6.4(1), $F$ must be two
dimensional.  Hence
$(1)$ is true.

Next, fix an extreme point $\xv$ of $I_s$.  If $\xv$ is not an extreme
point of $B$, then
it is properly contained in a minimal face $F$ of $B$. If $F$ were two
dimensional, then
by its minimality, we would get that $\xv$ lies in the relative interior of
$F$.  But in
this case $F\cap I_s$ would be a line segment by Lemma 6.4(2) and $\xv$
would be in its
interior, which is impossible since $\xv$ is extreme in $I_s$. Hence, $F$
has dimension
less than two.  Since $\xv$ is properly contained in $F$, this face must
have dimension
one.  Finally, since $\xv$ is not extreme in $B$, it must be contained  in
the relative
interior of $F$ and so $(2)$ is true.

Now suppose $\xv$ is an isolated extreme point in $I_s$.  Since $I_s$ has
dimension 2,
its faces consist of extreme points and line segments.  Since $\xv$ is an
isolated
extreme point,  it  must be the endpoint of adjacent line segments in the
boundary of
$I_s$. Hence $\xv$ lies in two distinct two dimensional faces of $B$ by
part $(1)$ of
this Lemma.  Thus part $(3)$ is true.

Finally, assume that the set $\e(B)$ of extreme points of $B$ is countable.
Observe
that each face in $B$ of dimension one is uniquely determined by its end
points, which are extreme points in $B$ by Theorem 0.2.  Hence $B$ has a
countable number of faces
of dimension one.

Now fix $s$ with $0\le s\le 1$.  We have that $\e(B)\cap I_s$ is countable.
Each
remaining extreme point of $I_s$ is contained in a unique face of $B$ with
dimension
one by part $(2)$ of this Lemma.   Since  there are a countable number of
one dimensional faces in $B$, the set of extreme points of $I_s$ is countable.

\end{proof}

\begin{lemma}   If $\xv$ is an isolated extreme point in $I_s$, then there are
central projections $z_1$ and $z_2$ in $B$ and $0 \le t \le 1$ such that
$z_1 \le z_2$ and
\[
\xv = \Psi(tz_1 + (1-t)z_2).
\]
\end{lemma}
\begin{proof}  As $\xv$ is isolated among the extreme points of $I_s$,
there are  distinct
two dimensional faces  $F_1$ and $F_2$ of $B$ such that $\xv \in F_1\cap
F_2 = F$ by
part $(3)$ of Lemma 6.5.  Next, since $F_1$ and $F_2$ have maximal
dimension in $B$, they
must be exposed and so we may write
\[
F_i = \Psi([\ptsmm{i},\ptspp{i}]), \quad i = 1,2.
\]
Further,  since the $F_i$'s intersect $I_s$ in adjoining line segments, we
must have that
$\tv_1 \ne \tv_2$.  Since neither of the $\tv_i$'s is zero, these vectors
are linearly
independent.  Hence, by Corollary 4.7, we must have
\[
F_3 = \Psi([z_1,z_2]),
\]
where $z_1$ and $z_2$ are central projections in $N$. Finally, since the
$F_i$'s have
dimension two, $F_3 = F_1\cap F_2$ must be a (possibly degenerate) line
segment.  Hence
elements of
$F_3$ are the images under $\Psi$ of convex combinations of $z_1$ and
$z_2$, as desired.

\end{proof}

With these preparations, we may now present the proof of Theorem 6.2.

\begin{proof} Write $\tv_1 = (1,0)$ and $\tv_2 = (0,1)$.  To show that
$N$ is abelian, it suffices to show that $p_{s,\tv_i}^\pm$ is central for
each spectral
pair of the form $(s,\tv_i),\,i = 1,2$.   Indeed in this case, by a
standard argument
in measure theory, we get that the spectral projections for $b_1$ and $b_2$
are central
and so $b_1$ and $b_2$ commute.

So now let $p$ denote such a projection, write $\xv = \Psi(p)$ and let
$I_s$ denote the isotrace slice containing $\xv$.  By part $(1)$ of Theorem
0.2, $\xv$ is an extreme
point of $B$ and so it is also an extreme point of $I_s$.  Since $\e(B)$ is
countable,
the extreme points of $I_s$ are countable by part $(4)$ of Lemma 6.5 and
since $I_s$ has
dimension 2, these extreme points form a  closed set.  Hence, the isolated
extreme points
of $I_s$ are dense in $I_s$ by Lemma 6.2.

Now let $\{\xv_n\}$ denote a sequence of isolated extreme point in $I_s$
that converges to
$\xv$.  Applying Lemma 6.6, we get sequences $\{z_{1,n}\}, \{z_{2,n}\}$ of
central
projections in $N$ and a sequence $\{t_n\}$ of real numbers in $[0,1]$ such
that
\[
\xv_n = \Psi(t_nz_{1,n} + (1-t_n)z_{2,n}),\quad n = 1,2,\dots\,.
\]
 Write $a_n = t_nz_{1,n}+(1-t_n)z_{2,n}$ for each $n$.  As the sequence
$\{a_n\}$ lies in
$N_1^+$, it is bounded and so it admits a subnet $\{a_\ga\}$ that converges
to an element
$a$ in $N_1^+$  in the weak$^*$-topology.  Since $\tau $ is normal, we have

\begin{align*}
\lim_\ga \Psi(a_\ga) &= \lim_\ga(\tau(a_\ga),\tau(b_1a_\ga),\tau(b_2a_\ga))\\
& = (\tau(a),\tau(b_1a), \tau(b_2a)) =  \Psi(a)\\
& =\xv = \Psi(p).
\end{align*}
Since $\Psi^{-1}(\xv) = \{p\}$ by part $(1)$ of Theorem 0.2, we get $a =
p$.   Since each $a_n$ is central and the center of $N$ is closed in the
weak$^*$--topology, we get that $p$ is central, as desired.

\end{proof}


\begin{thebibliography}{9}
\bibitem{AAW}  Charles A.\ Akemann, Joel\ Anderson and Nik\ Weaver,
{\em A Geometric Spectral Theory for n-tuples of Self-Adjoint Operators in
Finite
von Neumann Algebras},  { J.\ Functional\ Analysis
\textbf 164} (1999), 258-292

\bibitem{Dix} Jacques\ Dixmier, {\em $C^*-$ Algebra}, North-Holland, New York,
1977.


\bibitem{Rock} R. Tyrrell\ Rockafellar, {\em Convex Analysis}, Princeton
Univ. Press,
Princeton, 1970.

\bibitem{Rud:a}Walter\ Rudin, {\em Functional Analysis}, 2$^d$ Edition,
McGraw-Hill, New York, 1991.


\bibitem{Rud:b} Walter\ Rudin, {\em Real and Complex Analysis},
McGraw-Hill, New York,
1987.


\bibitem{Web} Roger\ Webster, {\em Convexity}, Oxford University Press, Oxford,
1994.
\end{thebibliography}
\end{document}